\documentclass[conference]{IEEEtran}
\ifCLASSINFOpdf
\else
\fi
\usepackage{multirow}
\usepackage{latexsym}
\usepackage{amsmath}
\usepackage{amsfonts}
\usepackage{amssymb}
\usepackage{graphicx}
\usepackage{times}
\usepackage{cite}

\begin{document}
%
\title{Pricing Energy in the Presence of Renewables}

\author{\IEEEauthorblockN{Ashkan Zeinalzadeh}
\IEEEauthorblockA{Department of Electrical Engineering\\
University of Notre Dame\\
azeinalz@nd.edu
}
\and
\IEEEauthorblockN{Indraneel Chakraborty}
\IEEEauthorblockA{Department of Finance\\
School of Business Admin\\
University of Miami\\
ichakraborty@bus.miami.edu
}
\and
\IEEEauthorblockN{Vijay Gupta}
\IEEEauthorblockA{Department of Electrical Engineering\\
University of Notre Dame\\
vgupta2@nd.edu
}}

\maketitle

\begin{abstract}
At present, electricity markets largely ignore the fact that renewable power producers impose significant externalities on non-renewable energy producers. This is because consumers are generally guaranteed electricity within certain load parameters. The intermittent nature of production by renewable energy producers implies that they rely on non-renewable producers so that the aggregate power delivered meets the promised quality of service. This implicit insurance provided by the non-renewable power sector to consumers is not currently priced and leads to an often ignored, hidden monetary transfer from non-renewable producers to renewable producers. As the fraction of energy supplied by renewable resources increases, these externalities also increase. In this paper, we quantify these externalities by developing the market clearing price of energy in the presence of renewable energy. We consider a day-ahead electricity market where renewable and non-renewable generators bid by proposing their asking price per unit of energy to an independent system operator (ISO). The ISO's problem is a multi-stage stochastic optimization problem to dispatch energy from each generator to minimize the cost of purchased energy on behalf of the consumers. We incorporate the notion of load variance using the Conditional Value-at-Risk (CVAR) measure in the day-ahead electricity market to ensure that the generators are able to meet the load within a desired confidence level. We analytically derive the market clearing price of energy as a function of CVAR. It is shown that a higher penetration level of the renewable energies may increase the market clearing price of energy.
\end{abstract}

\IEEEpeerreviewmaketitle

\section{Introduction}
\begin{figure}[!htb]
  \centering
  \includegraphics[width=68mm]{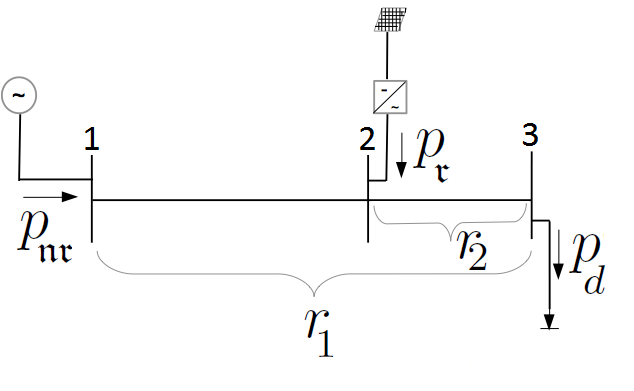}\\
  \caption{Electric grid model; $r_1$ is the resistance from bus~$1$ to the bus $3$; $r_2$ is the resistance from bus~$2$ to the bus $3$; $p_{\mathfrak{nr}}$ is the non-renewable power; $p_{\mathfrak{r}}$ is the renewable power; $p_{d}$ is the load power. }\label{fig:fig1}
\end{figure}

We consider an electric grid that delivers demanded electricity to consumers. The grid consists of a non-renewable generator unit, a renewable generator unit, a transmission line, and a consumer. The non-renewable generator unit is composed of $N$ generators with different output power limits, start-up costs, no-load costs, fuel costs per unit of generated energy, and minimum up and down times. The non-renewable and renewable generators are operated by different agents. For simplicity, we assume that there is only one transmission line and that the consumers are not able to distinguish between renewable energy (RE) and non-renewable energy (NRE). The renewable energy production is described by a stochastic process that may implicitly include curtailment. We assume that the grid takes all resultant renewable energy. Thus the non-renewable generator is dispatched to meet the net-load, which is the consumer load minus the renewable generator output (Figure~\ref{fig:fig1}).

The day-ahead energy market is designed to commit the non-renewable generators $24$-hours in advance and to set prices for each hour.
The load must be met within a confidence level.
The day-ahead market is cleared so as to minimize the total cost of energy. In the day-ahead bid, the non-renewable generators provide information about the asking price per unit of energy for each generator and their output limits to an independent system operator (ISO). The renewable generator only provides its demanded price for each unit of energy to the ISO.
The ISO after receiving the price of non-renewable and renewable energies, determines the output power for
each non-renewable generator that will ensure a given level of reliability to meet the load, and the corresponding market-clearing prices through the next $24$ hours.
The main result is that the market clearing price of energy may increase with greater penetration of renewable energy into the grid. Consequently, if the consumers impose the constraint of a maximum price for energy that they are willing to pay, a monetary transfer from the renewable producer to the non-renewable producer may be required to make the non-renewable producer whole. This transfer characterises the externality that the renewable producer imposes on the grid and is currently absorbed by the non-renewable producer. Further, if the price paid by the consumer is especially low, this may be a natural cap on the level of renewable penetration. The market clearing price is shown to be a non-decreasing function of uncertainty of the net load, the reliability demanded by the consumer, and the loss in the grid. Although the model we consider is stylized, the conclusions are striking and may point to a fundamental rethinking of the way renewable portfolio standards are currently planned.

Several works that are representative of the direction of this study include \cite{Bitar}-\cite{Joskow}, which study the problem under different market settings. The studies most related to our work are \cite{Bitar} and \cite{Bitar1}. They develop optimal strategies to inject wind energy into the grid under a fixed market price of energy. Unlike these works, we develop the market clearing price of energy to quantify the effect of uncertainty of the load and renewable energies on the market clearing price of energy.

The rest of the paper is structured as follows. In Section II, we formulate and solve the optimization problem of the ISO with the goal of minimizing the cost of energy while ensuring the planned generators in the day-ahead market are able to meet the load within a desired confidence level. In Section III we show the numerical results obtained through simulations. Finally, concluding remarks are provided in Section IV.

\section{ISO problem}

\subsection{Notation}

The following notations are used in this work.

\begin{itemize}
\item $p_{\ell}^t:$ Total active power loss at time $t$
\item $p_{d}^t:$ Active power requested by consumer
\item $p_i^t:$ Output power of the $i$th non-renewable generator at time $t$
\item $p_{\mathfrak{nr}}^t:$ Total active power generated by non-renewable generator unit
\item $p_{\mathfrak{r}}^t:$ Active power generated by renewable source
\item $r_i:$ Transmission line resistance
\item $(.)^+=\max(.,0)$
\item $s^t:=p_d^t+r_2 (p_r^t)^2-p_r^t$
\item $\pi_i:$ Asking price per unit of energy of the $i$th generating unit
\item $\pi_{\mathfrak{r}}:$ Asking price per unit of renewable energy.
\end{itemize}

\subsection{Load and Renewable Energy}

It is assumed that the load and renewable resource output are dependent, random variables with known distributions.
Let $\mathfrak{T}=\{1,2,...,T\}$ be the index set. Define $(P_D, P_{\mathfrak{r}})$ as the random processes on the probability
space $( \Omega,\mathfrak{F}, \mathbb{P})$. $P_D=\{ P_D^t,  t \in \mathfrak{T} \}$ and $P_{\mathfrak{r}}=\{ P_{\mathfrak{r}}^t, t \in \mathfrak{T} \}$ represent the load and renewable resource power respectively. For a fixed $t$, and for all $\omega \in \Omega$, $P_D^t(\omega)$ and $P_{\mathfrak{r}}^t (\omega)$ are non-negative random variables with known continuous probability density functions. For the given $\omega \in \Omega$, $P_D^t$ and $P_{\mathfrak{r}}^t$ are deterministic functions of $t$ that denote the load realization and renewable resource power realization respectively, at time $t$, donated by $p_{d}^t$, and $p_{\mathfrak{r}}^t$.

\subsection{Power Flow Constraints}

We assume the non-renewable generator unit is composed of $N$ generators. Let $p_i^t$ be the output power of the $i$th generating unit at time $t$. Rapid changes in output power, which cause rapid changes in the generator temperature or physical design, may increase maintenance costs.
The output of non-renewable generators can be limited by the generator capacity, or constraints on the quantity of fuel and $CO_2$ emissions.
Each unit must obey an output limit such that for all $i=1,...,N$
\begin{align} \label{eq:p11}
p_i^{\min} \leq p_i^t \leq p_i^{\max}.
\end{align}
The total output of non-renewable generator units is given as
\begin{align}
p_{\mathfrak{nr}}^{t}=\sum_{i=1}^{N} p_i^{t}.
\end{align}
We assume that the voltage at each bus is equal to one. The total active power loss ($p_{\ell}^t$) is approximated as
\begin{equation}\label{eq:actpower}
p_{\ell}^t \approx r_1 (p_{\mathfrak{nr}}^{t})^2+ r_2 (p_{\mathfrak{r}}^{t})^2.
\end{equation}
The power balance equation is given as
\begin{equation}\label{eq:p9}
p_{\mathfrak{nr}}^t+p_{\mathfrak{r}}^t-p_{\ell}^t=p_{d}^t,
\end{equation}
with a desired level of reliability $\omega \in \Omega$.

\subsection{Day-Ahead Market}

Let $\pi_{\mathfrak{r}}$ be the asking price per unit of energy of the renewable generator unit. The price of output energy of renewable generators is defined based on the levelized cost of electricity (LCOE), \cite{LCOE}. LCOE is the cost of one unit of energy over the lifetime of the generator and its associated parts. Let $\pi_i$ be the asking price per unit of energy of the $i$th non-renewable generator. Without loss of generality, we assume that $ 0< \pi_{\mathfrak{r}} < \pi_1 < ... < \pi_N$.

Let $p_\mathfrak{nr}:=\{(p_1^t,...,p_N^t)\}_{t=1}^{24}$ and $p_\mathfrak{r}:=\{p_\mathfrak{r}^t\}_{t=1}^{24}$. The ISO decides on a $p_\mathfrak{nr}$ that minimizes the total cost of energy through the next
$24$ hours as follows
\begin{equation}\label{eq:power}
\min_{p_{\mathfrak{nr}}} \sum_{t=1}^{24} E_{P_{D},P_{\mathfrak{r}}} \Big{[}  \sum_{i=1}^{N} \pi_i p_i^t+ \pi_{\mathfrak{r}} p_{\mathfrak{r}}^{t} \Big{]},
\end{equation}
with respect to (\ref{eq:p11})-(\ref{eq:p9}). We assume that the ISO is not aware of the start up cost, no-load cost and minimum up and down time of the generators. ISO only receives the bidding price and power from non-renewable generators. Therefore, the multi-stage optimization problem (\ref{eq:power}) is disjointed through stages (time). Let $n^t=(p_{\ell}^t+p_{d}^t-p_{\mathfrak{nr}}^t-p_{\mathfrak{r}}^t)^{+}$ and $F_{n^t}$ be the cumulative distribution function of the $n^t$, we use the concept of Value-at-Risk (VaR) and Conditional Value-at-Risk (CVaR), \cite{Rockafellar}. $VaR_{\alpha}(n^t)$ determines the worst possible $n^t$ that may occur within a given confidence level $\alpha$. For a given $0 < \alpha < 1$, the amount of $n^t$ will not exceed $VaR_{\alpha}(n^t)$ with probability $\alpha$,
\begin{align}\label{eq:var}
VaR_{\alpha}(n^t)=\min \{z| F_{n^t}(z) \geq \alpha \}.
\end{align}
CVaR is defined as the conditional expectation of $n^t$ above the amount VaR$_\alpha$. Let $E$ denote the expectation over $n^t$.
\begin{align}\label{eq:ninteen}
CVaR_{\alpha}(n^t)=E [ n^t | n^t > VaR_{\alpha}(n^t) ],
\end{align}
\begin{align}\label{eq:twenty}
CVaR_{\alpha}(n^t)=\int_{-\infty}^{\infty}z dF_{n^t}^{\alpha}(z),
\end{align}
where
\[
    F_{n^t}^{\alpha}(z)=
\begin{cases}
    0,& \text{if }  z < VaR_{\alpha}(n^t) \\
    \frac{F_{n^t}(z)-\alpha}{1-\alpha}, & \text{otherwise}
\end{cases}.
\]
In the day-ahead market the objective is to plan for the generators such that they are capable of meeting the load within the confidence level. We write the condition (\ref{eq:p9}) as
\begin{align}\label{eq:twentyyya}
CVaR_{\alpha}(n^t)=0.
\end{align}
Let $p_{\mathfrak{nr}}^t=(p_1^t,...,p_N^t)$. We define $(\tilde{\mu}_i,\mu_i)$ and $\lambda^t$ as the lagrange multipliers corresponding to (\ref{eq:p11})
and (\ref{eq:twentyyya}). The lagrange function for the ISO problem is given as
\begin{align}\label{eq:lagrange}
L^t(p_{\mathfrak{nr}})=&E_{P_{D},P_{\mathfrak{r}}} \Big{[} \sum_{i=1}^{N}  \pi_{i} p_{i}^t+ \pi_{\mathfrak{r}} p_{\mathfrak{r}}^t \Big{]}+\lambda^t \Big{[} CVaR_{\alpha}(n^t) \Big{]} \\ \notag
&+\sum_{i=1}^{N} \mu_i \Big{[} p_i^t-p_i^{\max}  \Big{]}
+\sum_{i=1}^{N} \tilde{\mu}_i \Big{[} p_i^{\min}-p_i^t  \Big{]}.
\end{align}
\\
\textbf{Theorem~$1$:} Let $s^t:=p_d^t+r_2 (p_r^t)^2-p_r^t$, it is claimed that
\begin{equation}\label{eq:cvar0}
CVaR_{\alpha}(n^t)=r_1 (\sum_{i=1}^{N}p_i^t)^2-\sum_{i=1}^{N}p_{i}^t+CVaR_{\alpha}(s^t).
\end{equation}

Proof. The proof is omitted and will be provided in the journal version due to space constraints.

It is evident from (\ref{eq:cvar0}) that $CVaR_{\alpha}(n^t)$ is convex in $p_{\mathfrak{nr}}^t$, therefore the lagrange function
(\ref{eq:lagrange}) is convex in $p_{\mathfrak{nr}}^t$. By substituting (\ref{eq:cvar0}) in (\ref{eq:lagrange})
\begin{align}
L^t(p_{\mathfrak{nr}}^t)=
&E_{P_{D},P_{\mathfrak{r}}} \Big{[} \sum_{i=1}^{N}  \pi_{i} p_{i}^t+ \pi_{\mathfrak{r}} p_{\mathfrak{r}}^t \Big{]}\\ \notag
&+\lambda^t \Big{[} r_1 (\sum_{i=1}^{N}p_i^t)^2-\sum_{i=1}^{N}p_{i}^t+CVaR_{\alpha}(s^t) \Big{]}\\ \notag
&+\sum_{i=1}^{N} \mu_i \Big{[} p_i^t-p_i^{\max}  \Big{]} +\sum_{i=1}^{N} \tilde{\mu}_i \Big{[} p_i^{\min}-p_i^t  \Big{]}.
\end{align}
The necessary Karush-Kuhn-Tucker (KKT) conditions for the ISO's problem are
\begin{align}\label{eq:KKT1}
&\pi_i+2 \lambda^t r_1 \sum_{i=1}^{N} p_i^t-\lambda^t+\mu_i-\tilde{\mu}_i=0,\\
&r_1 (\sum_{i=1}^{N}p_i^t)^2+CVaR_{\alpha}(s^t)=\sum_{i=1}^{N}p_{i}^t,\\
&p_i^{\min} \leq p_i^t \leq p_i^{\max},\label{eq:rang} \\
&\mu_i (p_i^t-p_i^{\max})=0,\label{eq:ortho}\\
&\tilde{\mu}_i (p_i^{\min}-p_i^t)=0,\\
&\mu_i \geq 0, \tilde{\mu}_i \geq 0. \label{eq:KKT5}
\end{align}
In order to have a feasible solution for the ISO's problem ((\ref{eq:KKT1})-(\ref{eq:KKT5})), Assumption~$1$ is considered.\\
\\
\textbf{Assumption~$1$:} Let $p^{\min}=\underset {1 \leq k \leq N}{\min} \,\ p_k^{\min}$. It is assumed that
\begin{itemize}
\item[a)] $1-4 r_1 CVaR_{\alpha}(s^t) \geq 0$.
\item[b)] For all $t=1,...,T$
\begin{equation}
\begin{cases}
      p^{\min} \leq CVaR_{\alpha}(s^t) \leq \sum_{i=1}^{N} p_i^{\max}, & \text{if}\ r_1=0 \\
      p^{\min} \leq \frac{1 \pm \sqrt{1-4 r_1 CVaR_{\alpha}(s^t) }}{2 r_1}  \leq \sum_{i=1}^{N} p_i^{\max}, & \text{if}\ r_1>0
    \end{cases}.
\end{equation}
\end{itemize}

Let $p^t$ be the solution of $r_1 (p^t)^2-p^t+CVaR_{\alpha}(s^t)=0$. Because of part $a)$ of Assumption~$1$ the value of $p^t$ is real
and because of part $b)$ of Assumption~$1$, there exists an unique $1 \leq k \leq N$ such that
$\sum_{i=1}^{k-1} p_i^{\max}  < p^t \leq \sum_{i=1}^{k} p_i^{\max}$. The value of $p^t$ based on the value of $r_1$
is given as follows.

\begin{equation}\label{eq:power1}
\begin{cases}
p^t=\frac{1 \pm \sqrt{1-4 r_1 CVaR_{\alpha}(s^t) }}{2 r_1}, & \text{if} \,\ r_1>0   \\ 
p^t=CVaR_{\alpha}(s^t), & \text{if} \,\ r_1=0 
\end{cases}
\end{equation}

By solving (\ref{eq:KKT1})-(\ref{eq:KKT5}), the values of $\{\mu_i\}_{i=1}^{N}$, $\{\tilde{\mu}_i \}_{i=1}^{N}$
and $\lambda^t$ are given below based on the value of $p^t-\sum_{i=1}^{k-1} p_{i}^{\max}$.
\begin{itemize}
\item $p_k^{\min} \leq p^t-\sum_{i=1}^{k-1} p_{i}^{\max} \leq p_k^{\max}$
\begin{align}
&\mu_i >0, \tilde{\mu}_i=0, p_i^t=p_i^{\max} \text{ for  \,\ } i=1,...,k-1,\\
&\mu_k=0, \tilde{\mu}_k=0, p_k^t=p^t-\sum_{i=1}^{k-1} p_{i}^{\max}, \\
&\mu_i=0, \tilde{\mu}_i=0, p_i^t=0, \text{ for all } i=k+1,...,N,\label{eq:SOKKT5}\\
&\lambda^t=\frac{\pi_i+\mu_i-\tilde{\mu}_i}{1-2r_1 p^t}\\ \notag
& \,\ \,\ =\frac{\pi_i+\mu_i-\tilde{\mu}_i}{\sqrt{1-4 r_1 CVaR_{\alpha}(s^t) }}, \,\ \text{for all} \,\ i=1,...,k\\
&\lambda^t=\frac{\pi_k}{\sqrt{1-4 r_1 CVaR_{\alpha}(s^t) }}.  \label{eq:result0}
\end{align}
\end{itemize}

To develop the generators optimal output for the case that $ 0 < p^t-\sum_{i=1}^{k-1} p_{i}^{\max} < p_k^{\min}$, we consider Assumption~$2$.\\
\\
\textbf{Assumption~$2$:} It is assumed that for all $i=1,...,N$
\begin{equation}\label{eq:actpower}
p_i^{\min} < \min_{k} \{ p_k^{\max}-p_k^{\min}  \}.
\end{equation}

Assumption~$2$ ensures that if the first $(k-1)$th generators are operating at their maximum power and additional power is needed to meet the day ahead load, when it is less than $p_k^{\min}$, then generator $(k-1)$th
can lower its output power without violating its output power, such that the $k$th generator operates at its minimum power ($p_{k}^{\min}$).
This is proved in Lemma~$1$, shown below, and is drawn from Assumption~$1$.\\
\\
\textbf{Lemma~$1$:} If $0 < p^t-\sum_{i=1}^{k-1} p_{i}^{\max} < p_k^{\min}$ then
\begin{equation}
p_{k-1}^{\min} < p^t-\sum_{i=1}^{k-2}p_i^{\max}-p_{k}^{\min} < p_{k-1}^{\max}.
\end{equation}

Proof. The proof is omitted and will be provided in the journal version due to space constraints.

\begin{itemize}
\item $ 0 < p^t-\sum_{i=1}^{k-1} p_{i}^{\max} < p_k^{\min}$
\begin{align}
&\mu_i >0, \tilde{\mu}_i=0, p_i^t=p_i^{\max} \text{ for all \,\ } i=1,...,k-2,\\
& \mu_{k-1} =0, \tilde{\mu}_{k-1}=0, p_{k-1}^t=p^t-\sum_{i=1}^{k-2}p_i^{\max}-p_{k}^{\min}, \\
&\mu_k=0, \tilde{\mu}_k>0, p_k^t=p_k^{\min}, \\
&\mu_i=0, \tilde{\mu}_i=0, p_i^t=0, \text{ for all } i=k+1,...,N,\label{eq:SOKKT5}\\
&\lambda^t=\frac{\pi_i+\mu_i-\tilde{\mu}_i}{1-2r_1 p^t} \\ \notag
& \,\ \,\ =\frac{\pi_i+\mu_i-\tilde{\mu}_i}{\sqrt{1-4 r_1 CVaR_{\alpha}(s^t) }}, \,\ \text{for all} \,\ i=1,...,k\\
&\lambda^t=\frac{\pi_{k-1}}{\sqrt{1-4 r_1 CVaR_{\alpha}(s^t) }}.  \label{eq:result}
\end{align}
\end{itemize}
It is evident from (\ref{eq:result0}) and (\ref{eq:result}), that the market clearing price of energy ($\lambda^t$) is higher at the times that
$s^t$ has a heavier tail distribution. A heavier tail distribution leads to a higher value of $CVaR_{\alpha}(s^t)$ and larger index of $k$
in (\ref{eq:result0}) and (\ref{eq:result}). Similarly, a higher level of reliability (larger $\alpha$) leads to a higher market clearing price of energy.
The accuracy of market clearing price (\ref{eq:result0}) and (\ref{eq:result}) is heavily dependent on the accuracy of the tail distribution of $s^t$.
The tail distribution of $s^t$ depends on the load and renewable energy distributions and model of loss function. In the next section, more descriptive simulations are presented.

\section{Simulations} \label{SectionV}

\textbf{Setup:} We consider that the non-renewable generator is composed of $6$ units. The asking price per unit of energy and the maximum output power for each unit is given in Table~\ref{table:tb0000}. It is assumed that $p_{i}^{\min}=0$ for all $i=1,...,6$. The level of reliability $(\alpha)$ demanded by consumer is fixed at $0.9$.

\begin{table}[h!]
\caption{Maximum output power ($p_i^{\max}$) and asking price ($\pi_i$)}
\label{table:tb0000}
\begin{center}
\begin{tabular}{|c|c|c|c|c|c|c|}
\hline
  $p_i^{\max}$ &0.05&0.1&0.12&0.15&0.18&0.25\\
\hline
   $\pi_i$  &20 & 30 & 40 & 50 & 60 & 70\\
\hline
\end{tabular}
\end{center}
\end{table}

\begin{figure}[h!]
\includegraphics[width=9cm]{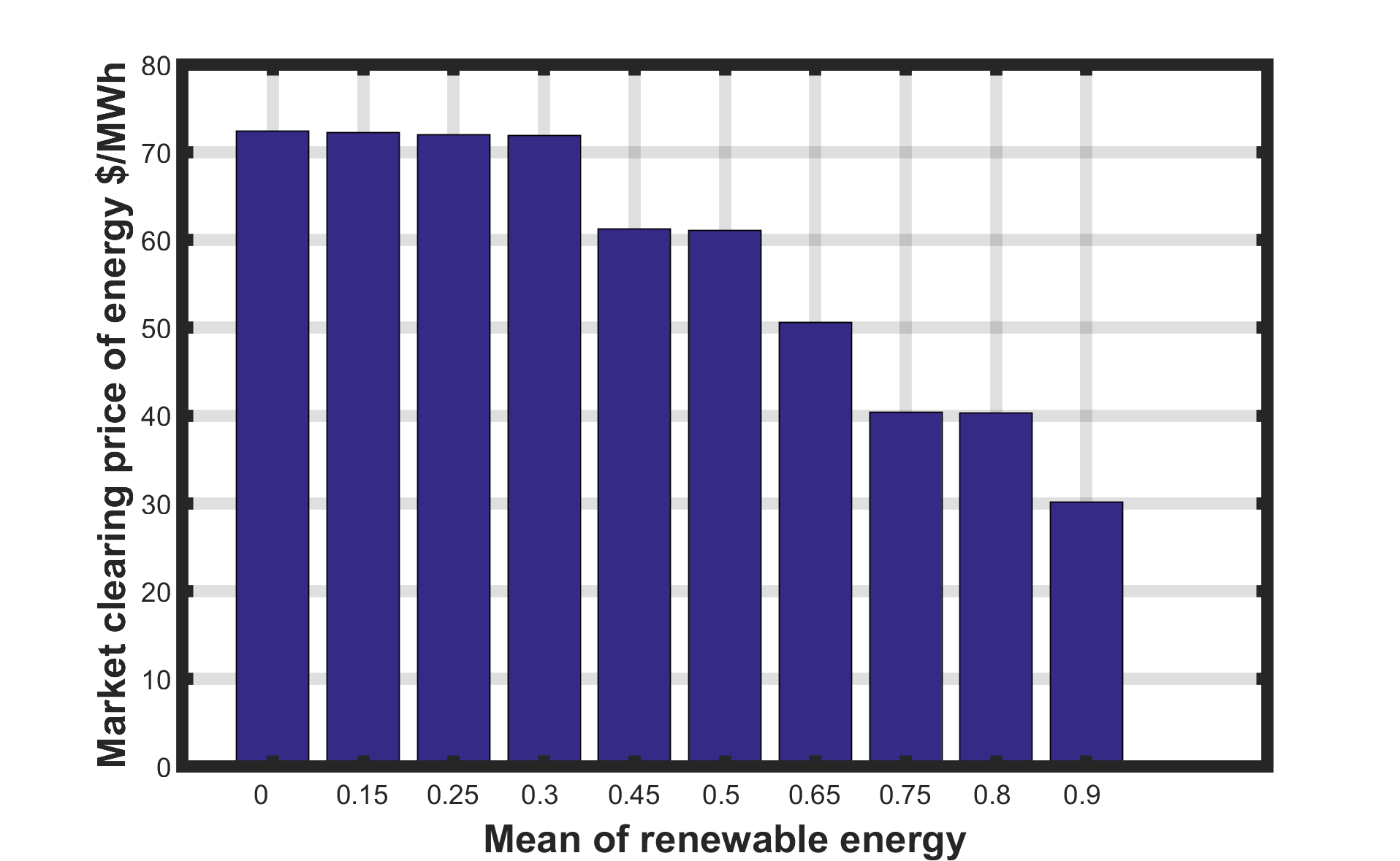}
\caption{Market clearing price versus the mean of renewable energy penetration.}
\label{fig:MEAN}
\end{figure}

We assume the load has a Gaussian distribution with a mean of $0.7$ and standard deviation of $0.1$. We repeat the simulation analysis for different scenarios of renewable energy penetration. In all scenarios, the renewable energy has a Gaussian distribution.\\
\\
\textbf{Case I:} We assume the standard deviation of the renewable energy is fixed at $0.1$ and the mean of renewable energy takes values $\{0, 0.15, 0.25, 0.3, 0.45, 0.5, 0.65, 0.75, 0.8, 0.9\}$. This corresponds to the naive expectation that the renewable energy penetration increases, while the uncertainty does not increase. The market clearing price of energy is plotted in Figure~\ref{fig:MEAN}.

\begin{figure}[h!]
\includegraphics[width=9cm]{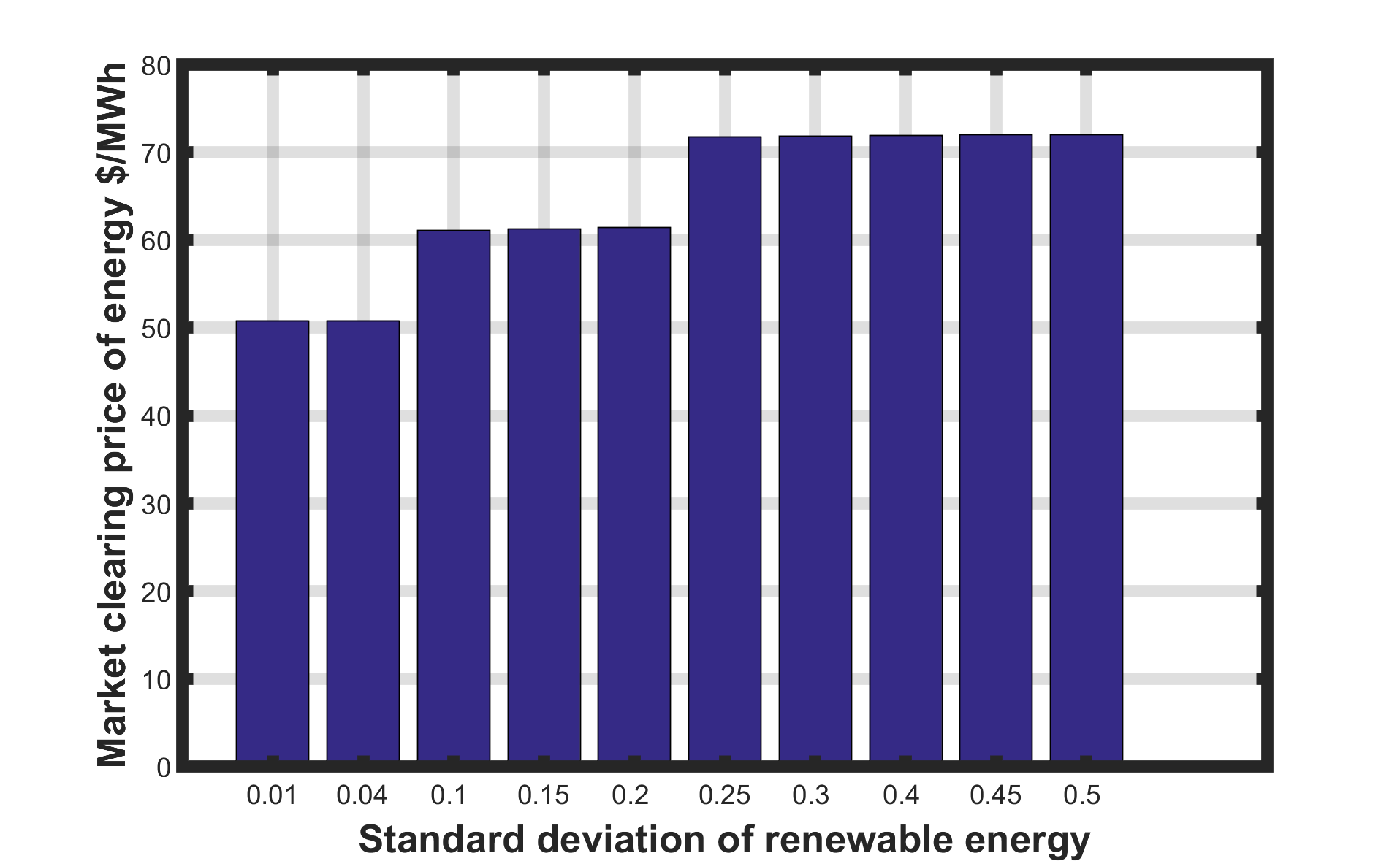}
\caption{Market clearing price versus the variance of renewable energy penetration.}
\label{fig:STD}
\end{figure}

It is observable from Figure~\ref{fig:MEAN} that the market clearing price of energy is decreasing in renewable energy penetration, if the standard deviation (i.e., the uncertainty) of renewable energy is kept constant while increasing the penetration. In reality of course, as more renewable energy is integrated the uncertainty in the total renewable energy also increases. We consider that scenario in Case III and IV below.\\
\\
\textbf{Case II:}
We assume that the mean of renewable energy is fixed at $0.5$ and the standard deviation of renewable energy takes values $\{0.01, 0.04, 0.1, 0.15,0.2,0.25,0.3,0.4,0.45,0.5\}$. The market clearing price of energy is plotted in Figure~\ref{fig:STD}.

\begin{figure}[h!]
\includegraphics[width=9cm]{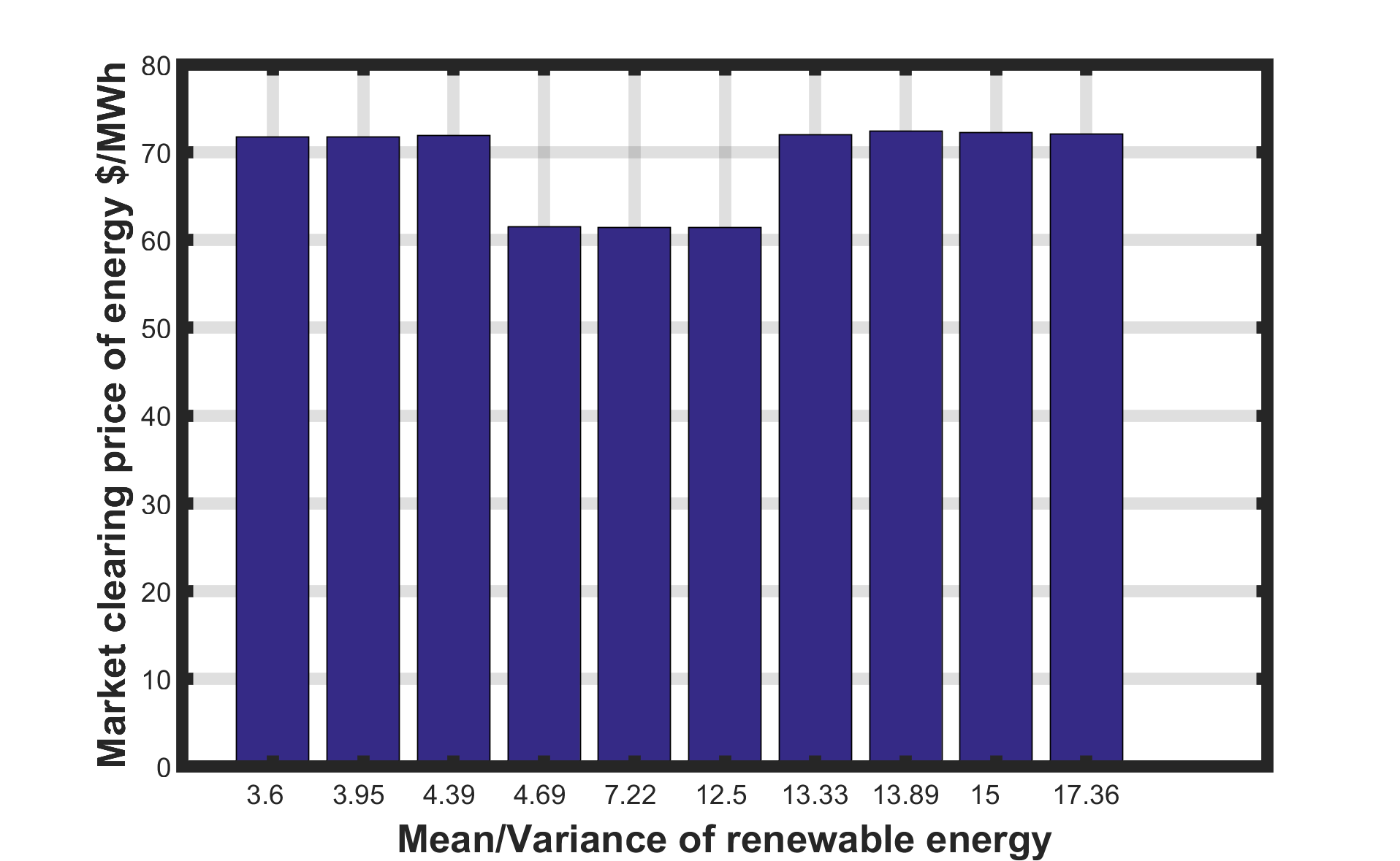}
\caption{Market clearing price versus the mean/variance of renewable energy penetration.}
\label{fig:STDMEAN}
\end{figure}

It can be seen from Fig.~\ref{fig:STD} that higher variability in renewable energy production increases the uncertainty in the net load. This increases the risk of unnecessary high capacity planning for the non-renewable generators in the day-ahead market, and leads to a higher market clearing price of energy.\\
\\
\textbf{Case III:}
Let the mean and standard deviation of renewable energy correspond to values \\
$\{0.05, 0.15, 0.25, 0.3, 0.45, 0.5, 0.65, 0.75, 0.8, 0.9\}$ and $\{0.06, 0.1, 0.12, 0.15, 0.32,0.2,0.3,0.4,0.45,0.5\}$ respectively. The market clearing price of energy is plotted in figure~\ref{fig:STDMEAN}.

\begin{figure}[h!]
\includegraphics[width=9cm]{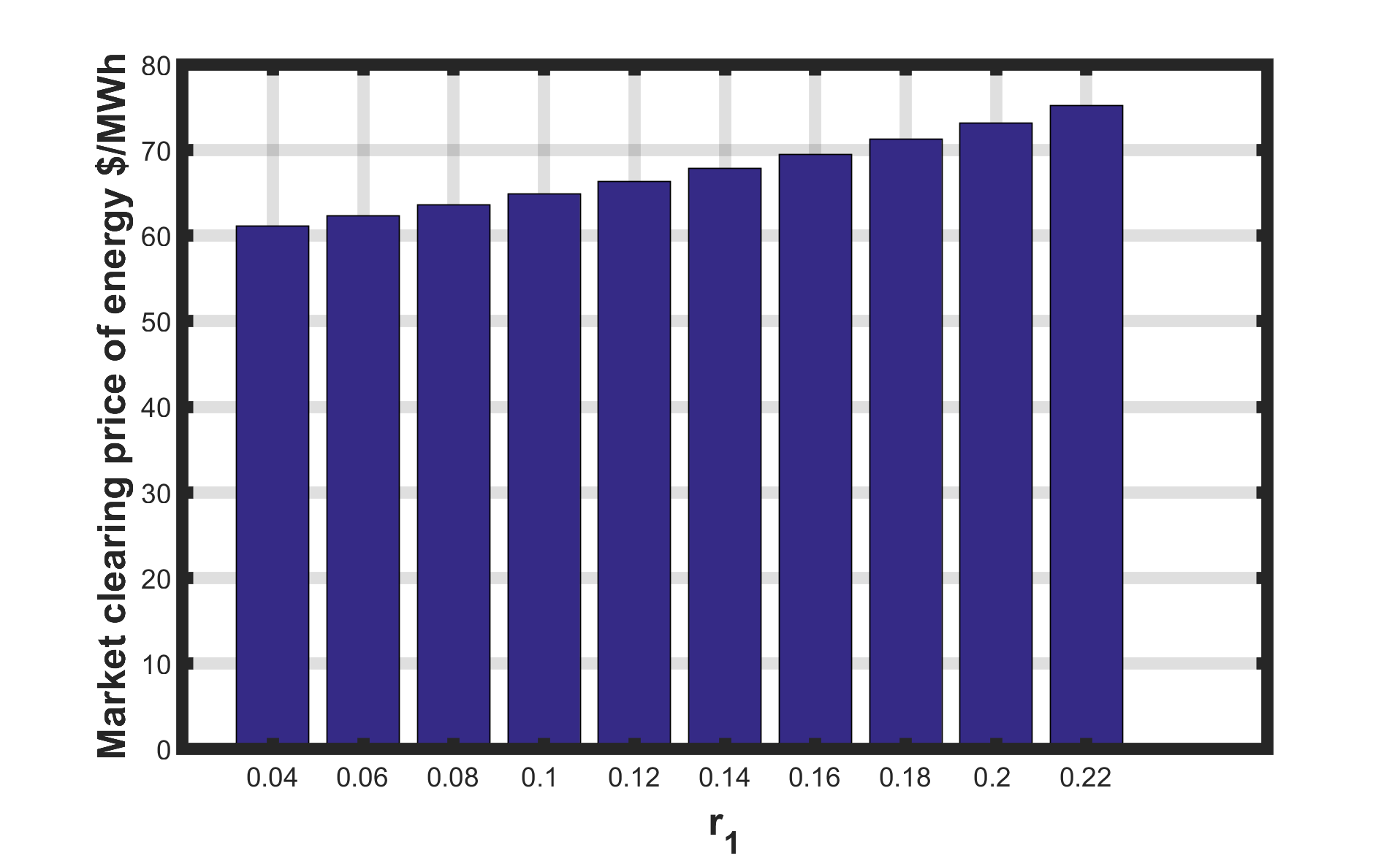}
\caption{Market clearing price versus the value of ($r1$).}
\label{fig:r1}
\end{figure}

In Figure~\ref{fig:STDMEAN}, the market clearing price of energy decreases until a certain level of the penetration is reached, after which the price increases. The price decreases in the beginning because of the lower marginal cost of renewable energy. However, after a certain level, the payment to the non-renewable generator to maintain the reliability constraints catches up and the market clearing price increases. The plot shows that (i) if the consumer insists on the lowest possible market clearing price, then the penetration level of renewable energy is capped; and (ii) if the consumer insists on a given price for the energy it may become important for the renewable producer to pay the non-renewable to compensate the latter.\\
%
%
%
\\
\textbf{Case IV:}
We assume the mean and standard deviation of renewable energy is fixed at $0.5$ and $0.1$ respectively. The line resistance $r_1$ takes values $\{0.04, 0.06, 0.08, 0.1, 0.12, 0.14, 0.16, 0.18, 0.2, 0.22\}$. The market clearing price of energy is plotted in figure~\ref{fig:r1}.

It is observable from Figure~\ref{fig:r1} that the market price of energy is increasing in the line resistance. A longer grid line (higher resistance) increases the loss of energy and increases the required capacity for non-renewable generators. Locating the renewable generators closer to the load leads to a lower market price of energy. This observation provides alternate methods to increase the renewable penetration.

\hfill


\section{Conclusion}
We quantify the market clearing price of energy in a day-ahead market as a function of the desired level of reliability. The expectation that increasing the penetration of the renewable energy reduces the market price of energy would be untrue if the uncertainty in the the renewable energy generation increases by the higher penetration level. If the consumer insists on paying no more than a certain price for energy, the renewable producer must transfer money to the non-renewable producer to compensate the latter. This transfer is currently ignored, which implies a hidden subsidy from non-renewable to the renewable producer. Understanding and resolving such frictions to fully consider the effects of uncertainty and fluctuations of renewable energies are central to fully realizing the benefits of renewable energy.

\section*{Acknowledgment}
The first author would like to thank the Center for Sustainable Energy at Notre Dame for partial support for this work.


\begin{thebibliography}{1}





\bibitem{Bitar}
Bitar, Eilyan Y, Ram Rajagopal, Pramod P Khargonekar, Kameshwar Poolla, and Pravin Varaiya, Bringing wind energy to market, IEEE Transactions on Power Systems, 27, 1225--1235, 2012.

\bibitem{Bitar1}
E. Bitar, A. Giani, R. Rajagopal, D. Varagnolo, P. Khargonekar, K. Poolla, P. Varaiya, Optimal contracts for wind power producers in electricity markets, 49th IEEE Conference on Decision and Control (CDC), 1919-1926, 2010.

\bibitem{Bathurst}
Bathurst, Graeme N, Jennie Weatherill, and Goran Strbac, Trading wind generation in short term energy markets, IEEE Transactions on Power Systems, 17, 782--789, 2002.

\bibitem{Safarinejadian}
Behrooz Safarinejadian and Masihollah Gharibzadeh and Mohsen Rakhshan, An optimized model of electricity price forecasting in the electricity market based on fuzzy timeseries, Systems Science and Control Engineering, vol. 2, no. 1, pp. 677--683, 2014.


\bibitem{Matevosyan}
Julija Matevosyan, and Lennart Soder, Minimization of imbalance cost trading wind power on the short-term power market, IEEE Transactions on Power Systems, 21, 1396--1404, 2006.

\bibitem{Morales}
Morales, Juan M, Antonio J Conejo, and Juan Pserez-Ruiz, Short-term trading for a wind power producer, IEEE Transactions on Power Systems, 25, 554--564, 2010.

\bibitem{Ashkan}
Ashkan Zeinalzadeh, and Vijay Gupta, Minimizing risk of load shedding and renewable energy curtailment in a microgrid with energy storage, Submitted to IEEE American Control Conference, 2017. 


\bibitem{Joskow}
Paul L. Joskow, Comparing the costs of intermittent and dispatchable electricity generating technologies, American Economic Review, vol 100, 238--241, 2011.

\bibitem{LCOE}
Projected costs of generating electricity, Nuclear Energy Agency, International Energy Agency, Organization for Economic Cooperation and Development, 2005.

\bibitem{Rockafellar}
Rockafellar, R. T. and S. Uryasev, Optimization of conditional value-at-risk. Journal of Risk, 2, pp. 21–41, 2000.

\end{thebibliography}
\end{document}